\documentclass[12pt]{article}

\usepackage{graphicx,tikz,color,url}
\usepackage{amsmath,amssymb,latexsym}
\usetikzlibrary{arrows}

\newtheorem{theorem}{Theorem}[section]

\newtheorem{lemma}[theorem]{Lemma}
\newtheorem{proposition}[theorem]{Proposition}
\newtheorem{observation}[theorem]{Observation}

\newtheorem{corollary}[theorem]{Corollary}
\newtheorem{remark}[theorem]{Remark}

\newtheorem{claim}{Claim}

\newcommand{\proof}{\noindent{\bf Proof.\ }}
\newcommand{\qed}{\hfill $\square$ \bigskip}
\newcommand{\cp}{\,\square\,}
\newcommand{\smallqed}{{\tiny ($\Box$)}}

\newcommand{\cH}{{\cal H}}

\makeatletter
\def\moverlay{\mathpalette\mov@rlay}
\def\mov@rlay#1#2{\leavevmode\vtop{%
   \baselineskip\z@skip \lineskiplimit-\maxdimen
   \ialign{\hfil$\m@th#1##$\hfil\cr#2\crcr}}}
\newcommand{\charfusion}[3][\mathord]{
    #1{\ifx#1\mathop\vphantom{#2}\fi
        \mathpalette\mov@rlay{#2\cr#3}
      }
    \ifx#1\mathop\expandafter\displaylimits\fi}
\makeatother

\newcommand{\cupdot}{\charfusion[\mathbin]{\cup}{\cdot}}

\begin{document}
	
	\title{Thresholds for the monochromatic clique transversal game}
	
	\author{Csilla Bujt\'as $^{a,b}$\thanks{Email: \texttt{bujtasc@gmail.com}} 
		\and Pakanun Dokyeesun $^{b}$\thanks{Email: \texttt{pakanun.dokyeesun@student.fmf.uni-lj.si}} 
		\and Sandi Klav\v zar $^{b,c,d}$\thanks{Email: \texttt{sandi.klavzar@fmf.uni-lj.si}}
	}
	\maketitle
	
	\begin{center}
		$^a$ Faculty of Information Technology, University of Pannonia, Hungary\\
		\medskip
		
		$^b$ Faculty of Mathematics and Physics, University of Ljubljana, Slovenia\\
		\medskip

		$^c$ Institute of Mathematics, Physics and Mechanics, Ljubljana, Slovenia\\
		\medskip
		
		$^d$ Faculty of Natural Sciences and Mathematics, University of Maribor, Slovenia\\
		\medskip
	\end{center}
	\begin{abstract}
		We study a recently introduced two-person combinatorial game, the $(a,b)$-monochromatic clique transversal game which is  played by Alice and Bob on a graph $G$.  As we observe, this game is equivalent to the $(b,a)$-biased Maker-Breaker game played on the clique-hypergraph of $G$. Our main results concern the threshold bias $a_1(G)$ that is the smallest integer $a$ such that Alice can win in the $(a,1)$-monochromatic clique transversal game on $G$ if she is the first to play.  Among other results,  we determine the possible values of $a_1(G)$ for the disjoint union of graphs, prove a formula for $a_1(G)$ if $G$ is triangle-free, and obtain the exact values of $a_1(C_n \cp C_m)$, $a_1(C_n \cp P_m)$, and $a_1(P_n \cp P_m)$ for all possible pairs $(n,m)$.
	\end{abstract}
	\noindent
	{\bf Keywords:} clique-hypergraph; Maker-Breaker game; clique transversal game,\\ threshold bias.\\
	
	
	\section{Introduction}
	
	The \emph{monochromatic transversal game} (MT game), as introduced by Mendes et al.\ \cite{mendes-2022+}, is a positional game on a hypergraph $\cH$. The two players alternately choose (play) an unplayed vertex and Alice colors it red while Bob colors the chosen vertex blue. Alice wins if she obtains a red transversal (vertex cover) in $\cH$ and Bob wins if he colors all vertices of a hyperedge with blue. It is called Alice-start game or Bob-start game, respectively, if Alice or Bob is the first to play. The \emph{monochromatic clique-transversal game} (MCT game) is played on the clique-hypergraph $\cH_G$ of a graph $G$. For the sake of simplicity, we will usually refer to an MCT game as being played on $G$ instead of $\cH_G$. The \emph{$(a,b)$-monochromatic clique-transversal game} ($(a,b)$-MCT game) follows the rules of the MCT game, but Alice (resp.\ Bob) colors $a$ vertices (resp.\ $b$ vertices) in each of her (resp.\ his) turn. 
	
	The introductory paper \cite{mendes-2022+} studied the $(a,b)$-MCT game on the clique-hyper\-graphs of powers of cycles and paths, observing also some basic properties of the game.\footnote{The recent conference paper \cite{mendes-2021} written by Mendes et al.\ started the study of $(a,b)$-MT game played on another type of hypergraphs, namely on the biclique hypergraphs of powers of paths and cycles.} In this paper, besides some general results we  concentrate  on triangle-free graphs in general and specifically on grids, cylinder, and torus graphs.  
	
	\subsection{Standard definitions}
	
	A \emph{hypergraph} $\cH$ is a set system over the vertex set $V(\cH)$. The hyperedge set $E(\cH)$ of $\cH$ may contain any nonempty subset of $V(\cH)$ i.e., $E(\cH) \subseteq 2^{V(\cH)} \setminus \{\emptyset\}$.  It is clear that every simple graph can be considered as a (2-uniform) hypergraph. A hypergraph is \emph{simple}, if no hyperedge $e \in \cH$ is a subset of a different $e'\in E(\cH)$. The \emph{maximum vertex degree $\Delta( \cH)$} of $\cH$ is the maximum number of hyperedges incident to a vertex in $\cH$. A \emph{singleton} is a hyperedge of cardinality $1$.
	
	A set $S\subseteq V(\cH)$ is a \emph{transversal} (also called vertex cover or hitting set) in the hypergraph $\cH$, if $S$ contains at least one vertex from each hyperedge $e \in E(\cH)$. Further, $S$ is a \emph{minimal transversal} if it contains no transversal as a proper subset. The \emph{transversal hypergraph} $Tr(\cH)$ of $\cH$ is defined on the vertex set $V(\cH)$ and the hyperedges correspond to the minimal transversals of $\cH$. As it has been stated already in Berge's fundamental book \cite{berge}, every simple hypergraph $\cH$ satisfies $\cH=Tr(Tr(\cH))$.
	\medskip
	
	In a graph $G$, a \emph{clique} is a complete subgraph that is inclusion-wise maximal. 
	The \emph{clique-hypergraph} $\cH_G$ of $G$ is defined on the same vertex set as $G$ while its hyperedges correspond to the vertex sets of the cliques in $G$, cf.~\cite{campos-2013, liang-2020}. Note that every  isolated vertex in $G$ corresponds to a singleton in $\cH_G$. Moreover, if $G$ is triangle-free and contains no isolated vertices, then $\cH_G \cong G$.
	
	For a non-negative integer $k$, a \emph{$k$-independent set} $S \subseteq V(G)$ is a  set of vertices in a graph $G$ such that the maximum degree of the subgraph induced by $S$ is at most $k$ \cite{caro-2013, chellali-survey, mao-2018}. The \emph{$k$-independence number} of $G$ is the maximum cardinality of a $k$-independent set in $G$, and is denoted by $\alpha_k(G)$. The $0$-independence number $\alpha_0(G)$ is the independence number $\alpha(G)$ of the graph. 
	
	The \emph{Cartesian product} $G\cp H$ of graphs $G$ and $H$ is defined on the vertex set $V(G)\times V(H)$ such that two vertices $(g,h)$ and $(g',h')$ are adjacent if either $gg'\in E(G)$ and $h=h'$, or $g=g'$ and $hh'\in E(H)$. If $h\in V(H)$, then the subgraph of $G\cp H$ induced by the vertex set $\{(g,h):\ g\in V(G)\}$ is isomorphic to $G$, called a {\em $G$-layer}, and denoted with  $G^h$. Analogously the $H$-layers are defined and denoted with $^g\!H$ for a fixed vertex $g\in V(G)$.  (See \cite{hik-2011} for more details on Cartesian products.) 
	\medskip
	
The {\em disjoint union} $G \cupdot H$ of (disjoint) graphs $G$ and $H$ has the vertex set $V(G)\cup V(H)$ and the edge set $E(G)\cup E(H)$. The {\em domination number} $\gamma(G)$ of a graph $G$ is the cardinality of a smallest dominating set, that is, a set $X\subseteq V(G)$ such that each vertex from $V(G)\setminus X$ has a neighbor in $X$. $P_n$ and $C_n$ respectively stand for the path and cycle graph on $n$ vertices. The order of the graph $G$ will be denoted by $n(G)$. The open neighborhood and the closed neighborhood of a vertex $v$ in $G$ will be respectively denoted by $N(v)$ and $N[v]$. For a set $S\subseteq V(G)$ we define $N[S] = \bigcup_{v\in S} N[v]$. For every positive integer $k$, $[k]$ stands for the set $\{1,\ldots,k\}$.  
	
	\subsection{Games}
	
	The \emph{Maker-Breaker game,} introduced by Erd\H os and Selfridge \cite{erdos-1973}, is played by Maker and Breaker on a hypergraph $\cH$. (For different instances of the Maker-Breaker game see~\cite{duchene-2020, glazik-2022, gledel-2020, stojakovic-2021}.) The vertex set $V(\cH)$ is called the \emph{board} and the hyperedges of $\cH$ are the  \textit{winning sets}. The two players alternately choose (i.e., play) an unplayed vertex. We say that a player \emph{claims} the set $S \subseteq V(\cH)$ if he (or she) plays all vertices from $S$. It is a \emph{Maker-start game} (resp.\ \emph{Breaker-start game}) if Maker (resp.\ Breaker) is the first player.
		\emph{Maker wins} the game if he claims a winning set, while \emph{Breaker wins} if she can prevent Maker from doing this. The latter equivalently means that Breaker plays at least one vertex from each winning set, that is, she claims a  minimal transversal of $\cH$. The following observation was stated in \cite{bujtas-2022}:
	\begin{observation} \label{obs:switch}
		A Maker-Breaker game on $\cH$ with a player A as Maker and a player B as Breaker is the same as the Maker-Breaker game on $Tr(\cH)$ where the roles of A and B are switched.
	\end{observation}
   By comparing the definitions, it is easy to see that an MT game on $\cH$ can always be considered as a Maker-Breaker game on $Tr(\cH)$ where Alice is Maker and Bob is Breaker. By Observation~\ref{obs:switch}, we may infer the following statements that give us a useful approach in our study.
   \begin{observation} \label{obs:switch-MCT}
   \begin{itemize}
   	\item[$(i)$] An MT game on a hypergraph $\cH$ is the same as the Maker-Breaker game on $\cH$ where Maker corresponds to Bob and Breaker corresponds to Alice.
   		\item[$(ii)$] An $(a,b)$-MCT game on $G$ is the same as the $(b,a)$-biased Maker-Breaker game on $\cH_G$ where Maker is Bob and Breaker is Alice.
   	\end{itemize}
   \end{observation}

	For a graph $G$ and a fixed positive integer $\ell$, the \emph{threshold bias} (or shortly threshold) $a_\ell(G)$ (resp.\ $a_{\ell}'(G)$) denotes the smallest positive integer $a$ such that Alice can win in the Alice-start (resp.\ Bob-start) $(a,\ell)$-MCT game. By definition, if $a <a_\ell(G)$, then Bob has a winning strategy in the Alice-start $(a, \ell)$-MCT game on $G$. It is also known (see e.g., Exercise 3.6.2 in \cite{hefetz}) that for every $a \ge a_\ell(G)$, Alice wins in the Alice-start $(a, \ell)$-MCT game. The same is true for the threshold $a_\ell'(G)$ that is, Alice wins in the Bob-start $(a, \ell)$-MCT game on $G$ if and only if $a \ge a_\ell'(G)$.
	
	\paragraph{Structure of the paper.} In the next subsection, we continue by observing some further basic properties of the $(a,b)$-MCT game. Section~\ref{sec:triangle-free} is devoted to triangle-free graphs, where we prove general formulas for the thresholds $a_1(G)$ and $a_1'(G)$ if $G$ is triangle-free. Then, in Section~\ref{sec:disjoint-union}, we establish lower and upper bounds on $a_1(G_1 \cupdot G_2)$ in terms of $a_1(G_1)$ and $a_1(G_2)$ and afterwards, we show that all values between the lower and upper bounds are realizable. The theorems in Section~\ref{sec:cartesian_path_cycle} give a complete picture on the exact values of thresholds $a_1(C_n \cp C_m)$, $a_1(C_n \cp P_m)$, $a_1(P_n \cp P_m)$, and $a_1'(C_n \cp C_m)$, $a_1'(C_n \cp P_m)$,  $a_1'(P_n \cp P_m)$.

	\subsection{Preliminaries} \label{sec:prelim}
	
	It can be proved via different simple approaches   (see e.g., Proposition 2.1.6 in \cite{hefetz}) that if Maker (resp.\ Breaker) can win as a second player in a Maker-Breaker game on $\cH$, she (resp.\ he) can also win as a first player. It can be checked that the analogous proofs work for the $(a,b)$-biased Maker-Breaker game. By Observation~\ref{obs:switch-MCT}, we may state the same property for the $(a,b)$-MCT game that was already mentioned in \cite{mendes-2022+} as Remark 2.1.
    \begin{observation} \label{re:1}
        If there exists a winning strategy for Alice (resp.\ Bob) in the Bob-start (resp.\ Alice-start) $(a,b)$-MCT game on a graph $G$, then there exists a winning strategy for  Alice (resp.\ Bob) when she (resp.\ he) starts the game on $G$.
     	    \end{observation}
	Observation~\ref{re:1} and the definition of the threshold bias directly implies the following relation.
	\begin{observation} \label{obs:a_1-prime-not-prime}
	It holds for every graph $G$ that 
	\begin{equation}
		\label{eq:a_1-prime-not-prime}
		a_1'(G) \ge a_1(G)\,.
	\end{equation}
	\end{observation} 
	
	In the rest of the paper we will also refer to the following statements.	
		\begin{observation}\label{Breakerwin}
	\begin{itemize}
	\item[$(i)$] If a hypergraph $\cH$ contains no singletons, then Breaker has a winning strategy in the $(1,\Delta(\cH))$-biased Maker-Breaker game on $\cH$ no matter who is the first player.
	\item[$(ii)$] 	Let $G$ be an isolate-free graph. If each vertex belongs to at most $k$ cliques in $G$, then Alice has a winning strategy in the $(k,1)$-MCT game on $G$ no matter who is the first player. 
	\item[$(iii)$] If $G$ is an isolate-free graph and every vertex belongs to at most $k$ cliques in $G$, then $a_1(G)\le a_1'(G)\le k$.
	\end{itemize}
	\end{observation}
\noindent{\it Proof.\ }  In the $(1,\Delta(\cH))$-biased Maker-Breaker game on $\cH$, Breaker may reply to Maker's every move $v$ by playing an unplayed vertex from each hyperedge incident to $v$ (and some additional vertices if necessary). This strategy ensures that whenever Maker plays a second vertex from some hyperedge, at least one vertex from this hyperedge has been already played by Breaker. Then, as $\cH$ contains no singletons, Maker cannot win, no matter who starts the game. This proves $(i)$.

By Observation~\ref{obs:switch-MCT}(ii), we may conclude statement $(ii)$ from $(i)$. The last part $(iii)$ then follows by the definition of the thresholds. \qed

	\section{Thresholds for triangle-free graphs} \label{sec:triangle-free}
	In this section we focus on triangle-free graphs that do not contain isolated vertices. For such a graph $G$, the clique hypergraph $\cH_G$ is isomorphic to $G$ and therefore, by Observation~\ref{obs:switch-MCT}(ii), Bob wins the $(a,b)$-MCT game on $G$ if and only if he claims two adjacent vertices. 
	\begin{theorem} \label{thm:threshold}
		If\/ $G$ is a triangle-free graph that contains no isolated vertex, then 
		\begin{align*}
			(i)\quad & a_1(G) = \min_{k \geq 0} \max\{k, n(G)-\alpha_k(G)\} = \min_{X\subseteq V(G)} \max\{\Delta(G-X), |X|\}\,,  \\
			(ii)\quad & a_1'(G) =  \Delta(G)\,. 
		\end{align*}
	\end{theorem}
	\proof To establish~(i), we are going to prove three claims. 
	\begin{claim} \label{cl:1}
		If\/ $\displaystyle a \geq \min_{k \geq 0} \max\{k, n(G)-\alpha_k(G)\}$, then Alice wins in the $(a,1)$-MCT game.
	\end{claim}
	\noindent{\it Proof.\ }
	Alice chooses a maximum $k^*$-independent set $Y$ such that $$\max\{k^*, n(G)-\alpha_{k^*}(G)\}= \min_{k \geq 0} \max\{k, n(G)-\alpha_k(G)\}.$$
	Observe that, by our condition, the set $X=V(G)\setminus Y$ contains at most $a$ vertices. In the first turn, Alice plays all vertices from $X$ and further $a-|X|$ arbitrarily chosen vertices if $a >|X|$. Since $Y$ is a $k^*$-independent set, $\Delta(G[Y]) \le k^* \le a$; that is each vertex played by Bob in his first turn or later would have at most $a$ unplayed neighbors. Alice thus can reply to each move $v$ of Bob by claiming all the unplayed neighbors of $v$ (and some further vertices if $a$ is big enough). With this strategy, Alice ensures that Bob can never claim two adjacent vertices and Alice's moves form a vertex cover that is a clique transversal in the isolate-free graph $G$. \smallqed 
	
	\medskip 
	
	\begin{claim} \label{cl:2}
		If\/ $a < \displaystyle \min_{X\subseteq V(G)} \max\{\Delta(G-X), |X|\}$, then Bob wins in the $(a,1)$-MCT game.
	\end{claim}
	\noindent{\it Proof.\ } Let $A_1$ be the set of the $a$ vertices played by Alice in the first turn. Under the present condition, $ a < \max\{\Delta(G-A_1), |A_1|\}$. As $|A_1|=a$, we may infer $1 \le a < \Delta(G-A_1)$. Consequently, in his first turn, Bob can play a vertex $v$ such that $\deg_{G-A_1}(v) > a$. Then, Alice cannot claim all the unplayed neighbors of $v$ in her next turn and Bob wins in his second turn by claiming an unplayed neighbor of $v$. \smallqed
	\medskip 
	
	\begin{claim} \label{cl:3}
		For every graph $G$
		$$\min_{X\subseteq V(G)} \max\{\Delta(G-X), |X|\}= \min_{k \geq 0} \max\{k, n(G)-\alpha_k(G)\}.$$
	\end{claim}
	\noindent{\it Proof.\ } Let $k^*$ be the smallest integer such that 
	$$\max\{k^*, n(G)-\alpha_{k^*}(G)\}=\min_{k \geq 0} \max\{k, n(G)-\alpha_k(G)\}.$$
	We select a maximum $k^*$-independent set $Y^*$ and define $X^*=V(G)\setminus Y^*$. By the choice of $k^*$ and $Y^*$, we have $|X^*|= n(G)-\alpha_{k^*}(G)$ and $\Delta(G-X^*)= k^*$. This gives
	$$\max\{\Delta(G-X^*), |X^*|\} = \max\{k^*, n(G)-\alpha_{k^*}(G)\} = \min_{k \geq 0} \max\{k, n(G)-\alpha_k(G)\}$$
	and, as $X^* \subseteq V(G)$, we conclude 
	$$\min_{X\subseteq V(G)} \max\{\Delta(G-X), |X|\} \le \min_{k \geq 0} \max\{k, n(G)-\alpha_k(G)\}.$$
	\medskip
	
	To prove the other direction, let $X'$ be a set of vertices in $G$ such that
	\begin{equation} \label{eq:3}
		\max\{\Delta(G-X'), |X'|\}=\min_{X\subseteq V(G)} \max\{\Delta(G-X), |X|\}.\end{equation}
	From the possible candidates for the role of $X'$ we choose one of minimum cardinality. This ensures that $Y'=V(G)\setminus X'$ is a maximum $k'$-independent set with $k'=\Delta(G-X')$. Note that
	$|X'|=n(G)-\alpha_{k'}(G)$ also follows and (\ref{eq:3}) remains true. We therefore obtain 
	$$\min_{X\subseteq V(G)} \max\{\Delta(G-X), |X|\} = \max\{k', n(G)-\alpha_{k'}(G)\} \ge 
	\min_{k \geq 0} \max\{k, n(G)-\alpha_k(G)\}$$
	that completes the proof of the claim. \smallqed
	\medskip
	
	Claims \ref{cl:1}--\ref{cl:3} directly imply that the threshold stated  in~(i) holds for every triangle- and isolate-free graph.
	\medskip
	
	To prove~(ii), we consider the Bob-start $(a,1)$-MCT game on $G$. As the graph is isolate-free, Bob cannot win in his first turn. If $a < \Delta(G)$, Bob plays a vertex $v$ of maximum degree in the first turn. As $a<\deg(v)$, Alice cannot claim all vertices from $N(v)$ in her response and hence, Bob can win the game in his second turn by playing an unplayed vertex which is adjacent to $v$. This shows that $a_1'(G)\ge \Delta(G)$ holds.
	Under the present conditions, every clique of $G$ is an edge and therefore, Observation~\ref{Breakerwin}(iii) implies $a_1'(G) \le \Delta(G)$. This completes the proof of~(ii). 
	\qed
	
	Notice that, independently of Observation~\ref{re:1}, Theorem~\ref{thm:threshold} directly implies $a_1(G) \leq a_1'(G)$ for every isolate- and triangle-free graph $G$. Indeed, by selecting $X=\emptyset$ we obtain 
	$$a_1(G) \leq \max\{\Delta(G - \emptyset), |\emptyset|\}=\Delta(G)=a_1'(G).$$
	
	\begin{remark}
		If $G$ is triangle-free but contains $\ell \ge 1$ isolated vertices, then Bob wins the Bob-start $(a,1)$-MCT game by claiming a clique $K_1$ in his first turn. As it is true for every integer $a$, the threshold $a_1'(G)$ does not exist. If Alice wins the Alice-start game, she has to play all isolated vertices in her first turn. A further argument analogous to the proof of Theorem~\ref{thm:threshold} yields the following formula:
		\begin{align*}
			a_1(G) & = \min_{k \geq 0} \max\{k, \ell+ n(G')-\alpha_k(G')\},
		\end{align*}
		where $G'$ denotes the graph obtained from $G$ by deleting all isolated vertices. As $n(G)=n(G') + \ell$ and $\alpha_k(G)=\alpha_k(G')+ \ell$, we may conclude the following equality:
		\begin{align*}
		a_1(G) & = \min_{k \geq 0} \max\{k, \ell+ n(G)-\alpha_k(G)\}.
		\end{align*}
	\end{remark}

	\section{Disjoint union of graphs}
	\label{sec:disjoint-union}
	
	\begin{proposition}\label{prop:union}
		If  $G_1$ and $G_2$ are graphs, then
		\begin{equation} \label{eq:union}
			\max\{a_1(G_1), a_1(G_2)\} \le a_1(G_1\cupdot G_2) \le a_1(G_1)+ a_1(G_2)\,.
		\end{equation}
	\end{proposition}
	
	\proof 
	Let $a_1(G_1) = k_1$, $a_1(G_2) = k_2$, and set $G = G_1\cupdot G_2$. Assume without loss of generality that $k_1\ge k_2$. 
	
	Consider the $(a,1)$-MCT game on $G$. Assume first that $a<k_1$. In the first move Alice selects at most $a<k_1$ vertices from $G_1$,  hence Bob has the following winning strategy. He will play only on $G_1$ as long as it is possible, following his winning strategy on $G_1$. Since $a < a_1(G_1)=k_1$, Bob indeed has such a strategy. In this way, Bob also wins on $G$. It follows  that $a_1(G) \ge k_1$ which proves the lower bound. 
	
	To prove the upper bound, we need to show that Alice has a winning strategy if $a = k_1 + k_2$. In this case, Alice follows her winning strategies on $G_1$ and $G_2$. More precisely, in the first move, she selects $k_1$ vertices from $G_1$ and $k_2$ vertices from $G_2$ according to her strategies on $G_1$ and $G_2$, respectively. In the rest of the game, after each move of Bob in  $G_i$, $i\in [2]$, she replies by selecting optimally $k_i$ vertices in $G_i$ and arbitrary additional $k_{3-i}$ vertices from $G_i$. In this way Alice wins the game and hence $a_1(G) \le a_1(G_1)+ a_1(G_2)$. 
	\qed
	
	Recall that a {\em caterpillar} is a tree in which a single path is incident to (or contains) every edge. This single path of the caterpillar is called the {\em spine}~\cite{west-2001}. For $m\ge 1$ and $\ell\ge 0$ let $T_{m,\ell}$ be the caterpillar whose spine has $m$ vertices and to each vertex of the spine, exactly $\ell$ leaves are attached. 
	
	In the next result we determine the threshold $a_1$ for the caterpillars $T_{m,\ell}$. As $T_{m,0}$ is a path of order $m$ and $a_1(P_m)$ was already established in~\cite{mendes-2022+}, we will assume $\ell \ge 1$ here.
	
	\begin{theorem}\label{thm:cater}
		If $m$ and $\ell$ are positive integers, then 
		\begin{align*}
			a_1(T_{m,\ell}) = \begin{cases}
				m; &  m \le \ell,\\
				\ell; & \lfloor \frac{m}{2} \rfloor\le \ell < m,\\
				\ell +1; & \lfloor \frac{m}{3} \rfloor -1 \le \ell \le \lfloor \frac{m}{2} \rfloor -1,\\
				\ell +2; & \ell \le \lfloor \frac{m}{3} \rfloor -2.
			\end{cases} 
		\end{align*}
	\end{theorem}
	
	\proof 
	Let $X = \{v_1, \ldots, v_m\}$ be the set of the (consecutive) vertices in the spine of $T_{m,\ell}$. We distinguish four cases depending on the size of $\ell$ with respect to $m$. 
	
	\paragraph{Case 1:} $m \le \ell$.\\
	Note that $|X|=m$ and $\Delta(T_{m,\ell}-X)=0$, hence Theorem~\ref{thm:threshold}(i) implies that $a_1(T_{m,l})\le m$.
	
	To prove that $a_1(T_{m,\ell})\ge m$, we need to show that Bob wins the Alice-start $(a,1)$-MCT game on $T_{m,\ell}$, when $a < m$. Partition $V(T_{m,\ell})$ into sets $V_i$, $i\in [m]$, where $V_i$ contains $v_i$ and all the leaves attached to it. Let $A_1$ be the set of $a$ vertices played by Alice in her first move. Then there exists $i\in [m]$ such that $A_1 \cap V_i = \emptyset$. Let  Bob reply by playing the vertex $v_i$. Since $a < m\le \ell$, Alice cannot claim all the neighbors of $v_i$ in her next turn and hence Bob can win the game. We conclude that $a_1(T_{m,\ell})\ge m$.
	
	\paragraph{Case 2:} $\lfloor \frac{m}{2} \rfloor\le \ell < m$.\\
	To show that $a_1(T_{m,\ell})\le \ell$, let $Y = \{v_2, v_4, \ldots, v_{2 \lfloor \frac{m}{2}\rfloor}\} \cup Y'$, where $Y'$ contains $\ell - \lfloor \frac{m}{2}\rfloor$ arbitrary additional vertices of $T_{m,\ell}$. Then $|Y|=\ell$ and $\Delta(T_{m,\ell}-Y) =\ell$, hence using Theorem~\ref{thm:threshold}(i) we get $a_1(T_{m,l})\le \ell$.
	
	To prove $a_1(T_{m,\ell})\ge \ell$, we need to show that Bob wins the Alice-start $(a,1)$-MCT game on $T_{m,\ell}$, when $a < \ell$. For this sake partition $V(T_{m,\ell})$ just as in Case 1 into the sets $V_i$, $i\in [m]$. If $A_1$ is the set of $a$ vertices played by Alice in her first move, then since $a< \ell < m$ there exists $i\in [m]$ such that $A_1 \cap V_i = \emptyset$. Then Bob replies by playing $v_i$ and wins the game in the next turn. We conclude that $a_1(T_{m,\ell})\ge \ell$.
	
	\paragraph{Case 3:} $\lfloor \frac{m}{3} \rfloor -1 \le \ell \le \lfloor \frac{m}{2} \rfloor -1$.\\
	Since $\ell\ge 1$, in this case we have $m\ge 4$. 
	To show that $a_1(T_{m,\ell})\le \ell+1$, let $Z = \{v_3, v_6, \ldots, v_{3\lfloor \frac{m}{3}\rfloor}\}$.
	Then $|Z|\le \ell + 1$ and $\Delta(T_{m,\ell}-Z) \le \ell+1$. Therefore, $\max\{|Z|, \Delta(T_{m,\ell}-Z)\} \le \ell + 1$. Applying Theorem~\ref{thm:threshold}(i) again we obtain $a_1(T_{m,l})\le \ell + 1$.
	
	Next, consider the Alice-start $(a,1)$-MCT game on $T_{m,\ell}$, when $a < \ell+1$. Partition $V(T_{m,\ell})$ into sets $V_i$, $i\in [\lfloor\frac{m}{2}\rfloor]$, where $V_i$ contains $v_{2i-1}, v_{2i}$ and the leaves attached to both of them. When $m$ is odd,  $V_{\lfloor\frac{m}{2}\rfloor}$ also contains $v_{m}$ and all the leaves attached to it. Let $A_1$ be the $a$ vertices played by Alice in her first move. Since $a < \ell + 1 \le (\lfloor \frac{m}{2} \rfloor - 1) + 1 = \lfloor \frac{m}{2} \rfloor$, by the Pigeonhole principle there exists $i\in [\lfloor \frac{m}{2} \rfloor]$ such that $A_1 \cap V_i = \emptyset$. Then Bob can play the vertex $v_{2i-1}$. Since neither the vertex $v_{2i}$ nor the leaves attached to $v_{2i-1}$ have been played, Bob wins the game in his next move. We conclude that $a_1(T_{m,\ell})\ge \ell + 1$.	
	
	\paragraph{Case 4:} $\ell \le \lfloor \frac{m}{3} \rfloor -2$.\\
	Since $\ell\ge 1$, in this case we have $m\ge 9$. To prove the first direction of the statement, we refer to Theorem~\ref{thm:threshold} which implies
	$$a_1(T_{m,\ell}) \le a_1'(T_{m,\ell})=\Delta(T_{m,\ell})=\ell+2.$$
	
	
	Consider now the Alice-start $(a,1)$-MCT game on $T_{m,\ell}$, where $a < \ell+2$. Partition $V(T_{m,\ell})$ into sets $V_i$, $i\in [\lfloor\frac{m}{3}\rfloor]$, where $V_i$ contains $v_{3i-2},v_{3i-1}, v_{3i}$, and all their leaves. If $3 \nmid m$, then add to $V_{\lfloor\frac{m}{3}\rfloor}$ also the (one or two) vertices $v_j$, where $j > 3\lfloor\frac{m}{3}\rfloor$, as well as all the leaves attached to them. Let $A_1$ be the $a$ vertices played by Alice in her first move. Since $a <\ell+2 \le \lfloor \frac{m}{3} \rfloor$, there exists $i\in [\lfloor\frac{m}{3}\rfloor]$ such that $A_1 \cap V_i = \emptyset$. Then Bob can play $v_{3i-1}$. Since neither $v_{3i-2}$ nor  $v_{3i}$ nor any leaf attached to $v_{3i-1}$ were played by now, Bob wins the game in the next turn. We conclude that $a_1(T_{m,\ell})\ge \ell+2$.
	\qed
	
	Note that the first case of Theorem~\ref{thm:cater}, that is $a_1(T)=m$, also holds for caterpillars $T$ which are obtained from $T_{m,m}$ by attaching some further leaves (without any restrictions) to the vertices of the spine.
	
	\begin{proposition}
		\label{thm:realizations}
		If $1\le \ell \le k$ and $0\le i\le \ell$, then $a_1(T_{k,k+i} \cupdot T_{\ell,k+i}) = k+i$.
	\end{proposition}
	
	\proof
	For this proof set $T_1 = T_{k,k+i}$, $T_2 = T_{\ell,k+i}$, and $T=T_1 \cupdot T_2$. Let $v_1, \ldots, v_k$ and $u_1, \ldots, u_\ell$ be the consecutive vertices of the spines of $T_1$ and $T_2$, respectively.
	
	
	First we set 
	$$X=\{v_2, v_4, \ldots, v_{2\lfloor k /2\rfloor},u_2, u_4, \ldots, u_{2\lfloor \ell/2 \rfloor} \}$$ 
	and observe that $|X|= \lfloor k /2\rfloor + \lfloor \ell/2 \rfloor \le k$ and $\Delta(T-X)=k+i$ hold. Theorem~\ref{thm:threshold} then implies $a_1(T) \le \max \{k,k+i\}=k+i$.

	Consider now the Alice-start $(a,1)$-MCT game on $T$, when $a < k + i$. Let $A_1$ be the set of vertices played by Alice in her first move. Then $|A_1| = a$. Because $a < k + i \le k + \ell$, by an argument parallel to the one from the proof of Theorem~\ref{thm:cater} we infer  that there exists a non-leaf $x$ from $T$ (so $x=u_i$ or $x=v_j$) such that neither $x$ nor any leaves attached to $x$ belong to $A_1$. Bob can then play $x$, and Alice cannot claim all its neighbors in her second move. Thus Bob wins the game and $a_1(T) \ge k+i$ follows. 
	\qed

	The following result shows that both inequalities in (\ref{eq:union}) are sharp and, moreover, every integer between  the lower and the upper bound may be the value of $a_1(G_1 \cupdot G_2)$ for appropriately chosen graphs.
	\begin{corollary} \label{cor:realization}
		For every three integers $\ell$, $k$, and $p$ satisfying $1 \le \ell \le k\le p \le \ell+k$, there exist caterpillars $G_1$ and $G_2$ such that
		$$ \max\{a_1(G_1), a_1(G_2)\}=k, \quad a_1(G_1\cupdot G_2)=p, \enskip \mbox{and} \enskip a_1(G_1)+a_1(G_2)=k + \ell\, .$$
	\end{corollary}
	\proof Let us consider the caterpillars $G_1=T_{k,p}$ and $G_2=T_{\ell,p}$. As $1\le \ell \le k \le p$, $a_1(G_1)=k$ and $a_1(G_2)= \ell$ hold by Theorem~\ref{thm:cater}. This already shows $\max\{a_1(G_1), a_1(G_2)\}=k$ and $a_1(G_1)+a_1(G_2)=k + \ell$. To prove the second equality, we set $i=p-k$ and observe that $0 \le i \le \ell$ holds under the present conditions. It follows from Proposition~\ref{thm:realizations} that $a_1(G_1\cupdot G_2)=k+i=p$. This completes the proof.
	\qed
	
	\section{Cartesian product of paths and cycles}
	\label{sec:cartesian_path_cycle}

	Throughout this section, let $V(C_n) = [n]$ for $n\ge 3$, and let $V(P_n) = [n]$ for $n\ge 1$, where $1$ and $n$ are end vertices of $P_n$. Moreover, for $C_n$ we will consider indices modulo $n$. 
	
	\subsection{Torus graphs}
	
	In this subsection we prove: 	
	
	\begin{theorem} \label{thm:CnCm}
		If $m\ge n \ge 3$, then
		\begin{align*}
			a_1(C_n \cp C_m) = \begin{cases} 
				1; &  n = m = 3,\\ 
				2; & n = 3, m = 4,\\
				3; & n = 3, m \ge 5,\\
				4; & m\ge n\ge 4,
			\end{cases}
		\end{align*}
		and
		\begin{align*}
			a_1'(C_n \cp C_m) = \begin{cases} 
				1; & n = m = 3,\\ 
				3; & n = 3, m\ge 4,\\
				4; & m\ge n\ge 4.
			\end{cases}
		\end{align*}
	\end{theorem}

	\proof 
	Consider first the case $m\ge n \ge 4$ and recall that $V(C_n\cp C_m) = \{(i,j):\ i\in [n], j\in [m]\}$. Then $C_n\cp C_m$ is a 4-regular, triangle-free graph, hence Theorem~\ref{thm:threshold}(ii) implies $a'_1(C_n\cp C_m) =4$. Thus by~\eqref{eq:a_1-prime-not-prime}, $a_1(C_n\cp C_m) \le a'_1(C_n\cp C_m) =4$. 
	
	It remains to show that $a_1(C_n\cp C_m) \ge 4$. For this sake consider the following winning strategy for Bob in the Alice-start $(3,1)$-MCT game on $C_n\cp C_m$. Let $A_1$ be the set of three vertices played by Alice in her first move. Since we have assumed that $m\ge n \ge 4$, there exist $i\in [n]$ and $j\in[m]$ such that $A_1 \cap V((C_n)^j) = \emptyset$ and $A_1 \cap V(^i(C_m)) = \emptyset$. Then Bob plays $(i,j)$ in his first move. Since $(i,j)$ is adjacent to four unplayed vertices,  Alice cannot claim all these vertices in her next turn. Therefore Bob can win by playing an unplayed neighbor of $(i,j)$. This finishes the proof for $m\ge n \ge 4$. 
	
	In the rest of the proof we assume that $n=3$ and consider the $(a,1)$-MCT game on $C_3 \cp C_m$. 
	
	\paragraph{Case 1:} $m =3$.\\ 
	Consider the Bob-start $(1,1)$-MCT game on $C_3 \cp C_3$. By symmetry we may assume without loss of generality that the first move of Bob is $(1,1)$. Then Alice replies by the move $(2,1)$. Then, using symmetry again, we need to consider four cases for a possible second move of Bob: $(1,2)$, $(2,2)$, $(3,2)$, and $(3,1)$. If Bob plays $(1,2)$, Alice replies by $(1,3)$. Then the only layers in which Alice has not played yet are $(C_3)^2$ and $^3(C_3)$. Whatever Bob plays next, Alice replies in the layer $(C_3)^2$ and later she can finish the game with a move in $^3(C_3)$. In all the other possible second moves of Bob, Alice replies by playing $(1,2)$ and then the only layers in which Alice has not played yet are $(C_3)^3$ and $^3(C_3)$. Then it is straightforward to see that in all the cases Alice can win the game. We conclude that $a'_1(C_3 \cp C_3) = 1$ and by Observation~\ref{re:1} also $a_1(C_3 \cp C_3) = 1$. 
	
	\paragraph{Case 2:} $m =4$.\\
	Since each vertex is incident with three cliques (one $3$-cycle and two edges), 
	Observation \ref{Breakerwin}(iii) yields $a'_1(C_3\cp C_4) \le 3$. 
	
	Consider the Bob-start $(2,1)$-MCT game on $C_3 \cp C_4$. Let Bob choose $(2,2)$ in his first move. Then Alice must play the vertices $(2,1)$ and $(2,3)$, for otherwise Bob can claim an edge in the layer $^2(C_4)$ in his next move. Then Bob replies with $(1,2)$. As Alice cannot play all the three vertices $(3,2)$, $(1,1)$, and $(1,3)$ in her next move, Bob can claim a clique in his second move. Thus $a'_1(C_3\cp C_4) \ge 3$. We conclude that $a'_1(C_3\cp C_4) = 3$.
	
	Next we consider the Alice-start $(2,1)$-MCT game on $C_3 \cp C_4$. To show that $a_1(C_3\cp C_4) \le 2$, we will give a winning strategy for Alice. Let her start the game by playing  $(1,1)$ and $(2,3)$. Assume that Bob plays a vertex $(i,j)$. If $(i,j)$ is in $(C_3)^1$ or in $(C_3)^3$, then Alice responds by playing two neighbors of $(i,j)$ in $^i(C_4)$, and then Bob cannot claim any clique incident to $(i,j)$. As Alice has already played a vertex from each $C_3$, in the continuation, after any move $(i',j')$ of Bob, Alice can reply by choosing the (at most two) unplayed neighbors in $^{i'}(C_4)$, and she wins the game. In the second case suppose that Bob's first move $(i,j)$ is from $(C_3)^2$ or $(C_3)^4$. If $(i,j) \in \{ (1,2), (1,4), (2,2), (2,4)\}$, then Alice replies by choosing one neighbor of $(i,j)$ in $(C_3)^j$ and the neighbor of $(i,j)$ in $^i(C_4)$ which was not yet selected. One can now check easily that Alice can win the game. The remaining subcase is $(i,j) = (3,2)$ (or symmetrically $(i,j) = (3,4)$). Then Alice responds by choosing $(3,1)$ and $(3,3)$. In each case it is now easy to see that Alice will win the game. This shows that $a_1(C_3\cp C_4) \le 2$.
	
	It remains to prove that $a_1(C_3\cp C_4) \ge 2$ by providing a winning strategy for Bob in the $(1,1)$-MCT game. Assume that Alice first plays $(i,j)$. Then Bob replies by playing $(i',j')$, where $i'\ne i$ and $j'\ne j$, and can win the game in his second move. We conclude that $a_1(C_3\cp C_4) = 2$.
	
	\paragraph{Case 3:} $m \ge 5$.\\
	By Observation~\ref{Breakerwin}(iii), $a_1(C_3\cp C_m) \le 3$. 
	
	To show that $a_1(C_3\cp C_m) \ge 3$, we will show that Bob has a winning strategy in the Alice-start $(2,1)$-MCT game on $C_3\cp C_m$.
	Assume that Alice plays two vertices in her first move. Then there are two consecutive layers $(C_3)^{k}$ and $(C_3)^{k+1}$ such that no vertex in these layers was played by Alice. By symmetry we may assume that the vertices $(1,k-1)$ and $(2,k-1)$ were not selected by Alice. Then Bob replies by selecting the vertex $(1,k)$. This forces Alice to reply by playing $(1,k-1)$ and $(1,k+1)$. Then Bob plays $(2,k)$ which forces Alice to play $(2,k-1)$ and $(2,k+1)$. Now Bob wins by playing $(3,k)$. Hence $a_1(C_3\cp C_m) \ge 3$ and consequently $a_1(C_3\cp C_m) = 3$.
	
	By Observation~\ref{Breakerwin}(iii), we have $3 = a_1(C_3\cp C_m) \le  a_1'(C_3\cp C_m) \le 3$. We conclude that $a_1'(C_3\cp C_m) = 3$. 
	\qed

	\subsection{Cylinder graphs}
	
	In this subsection we consider cylinder graphs, that is, Cartesian products of paths by cycles. The result reads as follows. 
	
	\begin{theorem} \label{thm:CnPm}
		If $n \ge 3$ and $m\ge 2$, then
		\begin{align*}
			a_1(C_n \cp P_m) = \begin{cases} 
				2; & n = 3,\ 2\le m \le 5,\text{or\ }\\
				& n = 4,\ m=2,\\
				3; & n = 3,\ m \ge 6, \text{or\ } \\
				& n = 4,\ 3 \le m \le 5, \text{or\ } \\		  
				& n=5,\ 2 \le m \le 4,\text{or\ } \\
				& 6 \le n \le 9,\ m=2,3,\text{or\ } \\
				& n \ge 10,\ m=2,\\
				4; & \text{otherwise}.
			\end{cases}
		\end{align*}
		and
		\begin{align*}
			a_1'(C_n \cp P_m) = \begin{cases} 
				2; & n=3,\ m=2,\\
				3; & n=3,\ m \ge 3,\text{or\ }\\
				& n \ge 4,\ m=2,\\
				4; & n \ge 4,\ m \ge 3.
			\end{cases}
		\end{align*}
	\end{theorem}

	\begin{table}
		\begin{center}
			\begin{tabular}{ c  | c c c c c }
				$n\backslash m$ & \,2\, & \,3\, & \,4\, & \,5\, & $\ge 6$  \\
				\hline
				3 & 2 & 2 & 2 & 2 & 3   \\
				4 & 2 & 3 & 3 & 3 & 4   \\
				5 & 3 & 3 & 3 & 4 & 4   \\
				6 & 3 & 3 & 4 & 4 & 4   \\
				7 & 3 & 3 & 4 & 4 & 4   \\
				8 & 3 & 3 & 4 & 4 & 4  \\
				9 & 3 & 3 & 4 & 4 & 4  \\
				$\ge 10$ & 3 & 4 & 4 & 4 & 4
			\end{tabular}
\qquad\qquad
			\begin{tabular}{ c  | c c c c c }
				$n\backslash m$ & \,2\, & \,3\, & \,4\, & \,5\, & $\ge 6$  \\
				\hline
				3 & 2 & 3 & 3 & 3 & 3   \\
				4 & 3 & 4 & 4 & 4 & 4   \\
				5 & 3 & 4 & 4 & 4 & 4   \\
				6 & 3 & 4 & 4 & 4 & 4   \\
				7 & 3 & 4 & 4 & 4 & 4   \\
				8 & 3 & 4 & 4 & 4 & 4  \\
				9 & 3 & 4 & 4 & 4 & 4  \\
				$\ge 10$ & 3 & 4 & 4 & 4 & 4
			\end{tabular}

		\end{center}
		\caption{$a_1(C_n\cp P_m)$ (left) and $a_1'(C_n\cp P_m)$ (right)}\label{Table:CnPm}
	\end{table}
	
	The values for $a_1(C_n\cp P_m)$ and $a_1'(C_n\cp P_m)$ are summarized in Table~\ref{Table:CnPm}. 
	
	In the rest of the subsection we prove Theorem~\ref{thm:CnPm}. We begin with the following lemma that allows us to use the results from the previous subsection.

	\begin{lemma} \label{lem:P_m}
		If $n \ge 3$ and $m\ge 4$, then 
		$$a_1(C_n \cp P_m) \le a_1(C_n \cp C_m) \quad {\rm and} \quad a_1'(C_n \cp P_m) \le a_1'(C_n \cp C_m)\,.$$
	\end{lemma}
	\proof Let $P_m$ be obtained from $C_m$ after the removal of the edge $1m$. The clique hypergraph of $C_n \cp P_m$ is then obtained from the clique hypergraph of $C_n \cp C_m$ by removing the cliques 
	$\{ (i,m),(i,1)\}$ for every $i\in [n]$.
	
	Suppose now that Alice can win in the $(a,1)$-MCT game on $C_n \cp C_m$ and consider one of her optimal strategies. If she plays according to the same startegy in the $(a,1)$-MCT game on $C_n \cp P_m$, Bob cannot claim a  clique of $C_n \cp C_m$ and therefore, he cannot claim a  clique in $C_n \cp P_m$. It follows that Alice can also win the $(a,1)$-MCT game on $C_n \cp P_m$. This argumentation is valid no matter Alice or Bob starts the games and implies that both $a_1(C_n \cp P_m) \le a_1(C_n \cp C_m)$ and $a_1'(C_n \cp P_m) \le a_1'(C_n \cp C_m)$ are true.
	\qed
	
	\begin{proposition}\label{prop:C3Pm}
		If $m \ge 2$, then 
		\begin{align*}
			a_1(C_3 \cp P_m) = \begin{cases} 
				2; &  2 \le m \le 5,\\ 
				3; & {\rm otherwise},
			\end{cases}
		\end{align*}
		and
		\begin{align*}
			a'_1(C_3 \cp P_m) = \begin{cases} 
				2; &  m=2,\\ 
				3; & {\rm otherwise}\,. 
			\end{cases}
		\end{align*}
	\end{proposition}
	
	\proof 
	Assume that $m \ge 2$ and consider the $(a,1)$-MCT game on $C_3 \cp P_m$. We first suppose that Alice starts the game.
	
	\paragraph{Case 1:} $2 \le m \le 5$.\\
	To show that $a_1(C_3 \cp P_m) \ge 2$, we will provide a winning strategy for Bob in the Alice-start $(1,1)$-MCT game. Assume that Alice plays $(i,j)$. Since $m \ge 2$, there exists a vertex $(i',j')$ with $i' \neq i$ and $j'\ne j$. Then Bob replies by playing $(i',j')$ and he can win the game within his next two moves. We infer that $a_1(C_3 \cp P_m) \ge 2$.
	
	To show that $a_1(C_3 \cp P_m) \le 2$, we will provide a winning strategy for Alice in the Alice-start $(2,1)$-MCT game.  We distinguish several subcases depending on $m$.
	
	\medskip\noindent
	{\em Case 1.1:} $m = 2$.\\
	In $C_3\cp P_2$, each vertex is incident with exactly two cliques. Hence Observation~\ref{Breakerwin}(ii) yields the result. 
	
	\medskip\noindent
	{\em Case 1.2:} $m = 3$.\\
	Alice plays $(1,2)$ and $(2,1)$ in her first move. If Bob plays $(i,j)$ in $(C_3)^1$ or in $(C_3)^2$, then Alice replies by playing all unplayed neighbors of $(i,j)$ in $^i(P_3)$ (and a further arbitrary unplayed vertex if it is needed). If Bob plays $(i,3)$, where $i \in [3]$, then Alice plays one vertex in the layer $(C_3)^3$ and $(i,2)$ if it is possible. Otherwise, Alice can play an arbitrary unplayed vertex. Now, every vertex is incident to at most two cliques which have not been played by Alice. Therefore, if Alice plays optimally, Bob cannot claim any clique and Alice wins the game.

	\medskip\noindent
	{\em Case 1.3:} $m = 4$.\\
	By Theorem~\ref{thm:CnCm} and Lemma~\ref{lem:P_m}, we have $a_1(C_3 \cp P_4) \le a_1(C_3 \cp C_4)= 2$.
	
	\medskip\noindent
	{\em Case 1.4:} $m = 5$.\\
	Alice plays $(1,4)$ and $(2,2)$ in her first move. Assume that Bob replies with the vertex $(i,j)$. If $(i,j)\in X = \{(1,2), (2,4), (3,2), (3,3), (3,4)\}$, then Alice chooses two neighbors of $(i,j)$ in the layer $^i(P_m)$. When $(i,j)\notin X$ then, if it is possible,  Alice replies by choosing two  neighbors of $(i,j)$, one in $(C_3)^j$ and the other in $^i(P_m)$. Otherwise, Alice can play any unplayed vertex. If she continues playing according to this strategy, Alice prevents Bob from claiming a clique and she wins the game.
	
	We conclude that $a_1(C_3 \cp P_m) \le 2$ for $2 \le m \le 5$.
	
	\paragraph{Case 2:} $m \ge 6$.\\
	By Lemma~\ref{lem:P_m} and Theorem~\ref{thm:CnCm} we have $a_1(C_3 \cp P_m) \le a_1(C_3 \cp C_m)= 3$.
	
	To prove that $a_1(C_3 \cp P_m) \ge 3$ we show that Bob can win in the Alice-start $(2,1)$-MCT game. After Alice plays two vertices in her first move, there exist  three consecutive layers $(C_3)^{j-1}$, $(C_3)^{j}$, and $(C_3)^{j+1}$, such that at most one vertex from these three layers were played by Alice and the possibly played vertex is not from $(C_3)^{j}$. Let $(1,j-1)$ be the possibly played vertex. No matter whether it was played or not, Bob replies by $(2,j)$. Then Alice is forced to select $(2,j-1)$ and $(2,j+1)$. After that, Bob plays $(3,j)$ and will win in his next move.  
	
	\medskip	
	It remains to consider the Bob-start game on $C_3\cp P_m$. If $m=2$, each vertex is contained in two  cliques and Observation~\ref{Breakerwin}(iii)  implies $a'_1(C_3\cp P_2) \le 2$.  By~\eqref{eq:a_1-prime-not-prime} and the first part of this proof we have $a_1'(C_3\cp P_2) \ge a_1(C_3\cp P_2) = 2$. Therefore $a'_1(C_3\cp P_2) = 2$.
	
	Suppose that $m \ge 3$. Then $a'_1(C_3\cp P_m) \le 3$ by Observation~\ref{Breakerwin}(iii). To see that $a'_1(C_3\cp P_m) \ge 3$ consider the following strategy of Bob in the Bob-start $(2,1)$-MCT game. In his first move, he plays $(1,2)$. Then Alice is forced to play $(1,1)$ and $(1,3)$. Afterwards, Bob replies by the move $(2,2)$ and since Alice cannot play all of its neighbors, Bob can win in the next turn. 
	 
	\qed
	
	Next, we consider the thresholds for $C_n \cp P_2$ where $n \ge 4$. 
	\begin{proposition}\label{prop:CnP2}
		If $n \ge 4$, then 
		\begin{align*}
			a_1(C_n \cp P_2) = \begin{cases} 
				2; &  n=4,\\ 
				3; &  n\ge 5.
			\end{cases}
		\end{align*}
		\end{proposition}
	
	\proof Assume that $n \ge 4$ and consider the Alice-start $(a,1)$-MCT game on $C_n \cp P_2$.
	
	\paragraph{Case 1:} $n=4$.\\
	Let $X = \{(1,1), (3,2)\}$. Then $\Delta((C_4 \cp P_2) -X) = 2$. By Theorem \ref{thm:threshold}(i), we have that  $a_1(C_4 \cp P_2) \le \max\{\Delta((C_4 \cp P_2) -X), |X|\}=2$. It suffices to show that Bob has a winning strategy in the $(1,1)$-MCT game.
	Assume Alice plays $(i,j)$ in her first move. Then Bob replies by choosing a vertex $(i',j')$ such that $i' \neq i$ and $j' \neq j$, and he will win the game in his second move. Hence $a_1(C_4 \cp P_2)=2$.
	
	\paragraph{Case 2:} $n\ge 5$.\\
	Since each vertex is incident with three cliques,  Observation~\ref{Breakerwin}(iii) implies that $a_1(C_n \cp P_2) \le 3$. It remains to show that $a_1(C_n \cp P_2) \ge 3$ by giving a winning strategy for Bob in the $(2,1)$-MCT game.
	Let $A_1$ be the set of two vertices played by Alice in her first move.
	There exists a vertex $(i,j)$ such that $(i,j) \notin N[A_1]$. Then Bob plays $(i,j)$ in his first move and can win the game in his next turn. Thus $a_1(C_n \cp P_2) = 3$.
	\qed
	
For $n \ge 4$, the Cartesian product $C_n \cp P_m$ is triangle-free with no isolated vertices. By Theorem~\ref{thm:threshold}, $a_1'(C_n \cp P_m) = \Delta(C_n \cp P_m)$ which implies the following:

	\begin{observation}\label{obs:a'CnPm}
	If $n\ge 4$, then 
	\begin{align*}
		a_1'(C_n \cp P_m) = \begin{cases} 
			3; & m=2,\\
			4; & m\ge 3\,.
		\end{cases}
	\end{align*}		
\end{observation}

	\begin{lemma}\label{lem:CnPm3}
		Assume that $n \ge 4$ and $m \ge 3$. Then $a_1(C_n \cp P_m) = 3$ if and only if $\gamma(C_n \cp P_{m-2}) \le 3$.
	\end{lemma}	
	
	\proof Let $G$ be the induced subgraph of $C_n \cp P_m$ which remains after the deletion of the vertices from $V((C_n)^1) \cup V((C_n)^m)$. Then $G$ is isomorphic to $C_n \cp P_{m-2}$. It suffices to show that  $a_1(C_n \cp P_m) = 3$ if and only if $\gamma(G) \le 3$.
	
	$(\Leftarrow)$ Let $X$ be a minimum dominating set of $G$ and assume that $|X|\le 3$. We will show that $a_1(P_n \cp C_m) = 3$. Since each vertex in $V((C_n)^1)\cup V((C_n)^m)$ has degree $3$ in $C_n \cp P_m$, and $X$ is a dominating set of $G$, $\Delta ((C_n \cp P_m) - X) \le 3$. By Theorem~\ref{thm:threshold}(i), $a_1(C_n \cp P_m) \le \max\{\Delta ((P_n \cp P_n) - X), |X|\}\le 3$.
	To show that $a_1(C_n \cp P_m) \ge 3$, we will show that Bob has a winning strategy in the Alice-start $(2,1)$-MCT game. After Alice plays two vertices in her first move, there is a vertex $(i,j)$ such that $(i,j)$ is adjacent to at least three unplayed vertices. Then Bob plays $(i,j)$ and wins the game in his next move. Therefore $a_1(C_n \cp P_m) =3$.
	
	$(\Rightarrow)$ Assume that $a_1(C_n \cp P_m) = 3$. By Theorem~\ref{thm:threshold}(i), there exists a set $X$ such that $\Delta((C_n\cp P_m)) - X) \le 3$ and $|X|\le 3$.  Then all vertices of degree $4$ in $C_n \cp P_m$ belong to $N_{C_n \cp P_m}[X]$ which implies that $X$ dominates all vertices from $V(G)$. If there exists a vertex $(i,j) \in X-V(G)$, then $(i,j)$ is adjacent to exactly one vertex $(i',j')$ from $V(G)$ and we can replace $(i,j)$ with $(i',j')$ in $X$. By repeating this process for the remaining vertices of $X$ if necessary, we find a dominating set of $G$ of order at most $3$. 
	\qed
	
	
	In \cite{pavlic-2013}, the domination number of $\gamma(C_n \cp P_m)$ was determined for each pair of parameters with $n \ge 4$ and $ 2 \le m \le 7$. It is well-known that $\gamma(C_n) = \lceil\frac{n}{3}\rceil$   and that  $\gamma(C_n \cp P_m) > 3$ holds whenever $n \ge 4$ and $m \ge 8$.  From these facts we can deduce the following statement: 
	
	\begin{proposition}
		\label{prop:gamma-CnPm}	
		If $n\ge 4$ and $m\ge 1$, then $\gamma(C_n \cp P_m) \le 3$ if and only if one of the following conditions holds:
		\begin{itemize}
			\item[$(i)$] $n = 4$ and $m \in [3]$;
			\item[$(ii)$] $n=5$ and $m \in [2]$;
			\item[$(iii)$] $6\le n\le 9$ and $m=1$.
		\end{itemize} 
	\end{proposition}
	
We are now in a position to complete the proof of Theorem~\ref{thm:CnPm}. Since the other cases have already been covered, it sufficies to consider cylinder graphs $C_n\cp P_m$ with $n\ge 4$.

From Lemma~\ref{lem:CnPm3} and Proposition~\ref{prop:gamma-CnPm} we deduce that if $n\ge 4$, then $a_1(C_n\cp P_m) = 3$ if and only if either $n=4$ and $3\le m\le 5$, or $n=5$ and $3\le m\le 4$, or $6\le n\le 9$ and $m = 3$. It remains to show that $a_1(C_n \cp P_m) = 4$ in each of the following cases: 
\begin{itemize}
\item[$(i)$] $n = 4$ and $m \ge 6$, 
\item[$(ii)$] $n = 5$ and $m \ge 5$,
\item[$(iii)$] $n \ge 6$ and $m \ge 4$,
\item[$(iv)$] $n \ge 10$ and $m =3$.
\end{itemize}
For these cases, $a_1(C_n\cp P_m) \neq 3$ and Theorem~\ref{thm:threshold}(i) implies $a_1(C_n\cp P_m) >2$. From Observation~\ref{Breakerwin}(iii) we deduce $a_1(C_n\cp P_m) \le 4$. Therefore,  $a_1(C_n\cp P_m) = 4$ holds under the conditions $(i)$-$(iv)$. This finishes the proof of Theorem~\ref{thm:CnPm}.
	
	\subsection{Grid graphs}
	
		\begin{theorem} \label{thm:PnPm}
		If $m\ge n \ge 2$, then
		\begin{align*}
			a_1(P_n \cp P_m) = \begin{cases} 
				2; & n = 2, 2 \le m \le 5, \text{or\ } \\
				&  n=3, m=3,4,\\
				3; & n = 2, m \ge 6, \text{or\ } \\		  
				& n=3, 5 \le m \le 11,\text{or\ } \\
				& n=4, 4 \le m \le 7,\text{or\ } \\
				& n=m=5,\\
				4; & \text{otherwise}.
			\end{cases}
		\end{align*}
		and
				\begin{align*}
		a_1'(P_n \cp P_m) = \begin{cases} 
		2; & n=m=2,\\
		3; & n=2,m \ge 3,\\
		4; & m\ge n \ge 3\,.
	\end{cases}
				\end{align*}
	\end{theorem}
	
Table~\ref{Table:PnPm} summarizes the values of $a_1(P_n\cp P_m)$. 
		\begin{table}
		\begin{center}
			\begin{tabular}{ c  | c c c c c c c c c c c c }
				$n\backslash m$ & \,2\, & \,3\, & \,4\, & \,5\, &\,6\, &\,7\, &\,8\, &\,9\, &\,10\, &\,11\, & $\ge 12$  \\
				\hline
				2 & 2 & 2 & 2 & 2 & 3 & 3 & 3 & 3 & 3 & 3 & 3   \\
				3 &  & 2 & 2 & 3 & 3 & 3 & 3 & 3 & 3 & 3 & 4  \\
				4 &   &   & 3 & 3 & 3 & 3 & 4 & 4 & 4 & 4 & 4   \\
				5 &   &   &   & 3 & 4 & 4 & 4 & 4 & 4 & 4 & 4  \\
				$\ge 6$ &   &   &   &   & 4 & 4 & 4 & 4 & 4 & 4 & 4
			\end{tabular}
		\end{center}
		\caption{$a_1(P_n\cp P_m)$ for $2 \le n\le m$}\label{Table:PnPm}
	\end{table}
	
	In the rest of the subsection we prove Theorem~\ref{thm:PnPm}. Observe that the Cartesian product $P_n \cp P_m$ is triangle-free with no isolated vertices for all $n,m \ge 2$. By Theorem~\ref{thm:threshold}(ii), $a_1'(P_n \cp P_m) = \Delta(P_n \cp P_m)$ and thus, we can conclude that
	\begin{align*}
		a_1'(P_n \cp P_m) = \begin{cases} 
			2; & n=m=2,\\
			3; & n=2,m \ge 3,\\
			4; & m\ge n \ge 3\,.
		\end{cases}
	\end{align*}
	
	Next, we consider the threshold for the Alice-start game on $P_2 \cp P_m$  where $m \ge 2$.
	\begin{proposition}\label{prop:P2Pm}
		If $m \ge 2$, then 
		\begin{align*}
		a_1(P_2 \cp P_m) = \begin{cases} 
		2; & 2 \le m \le 5,\\
		3; & m\ge 6\,.
		\end{cases}
		\end{align*}
	\end{proposition}
	\proof Assume that $m \ge 2$ and consider the Alice-start $(a,1)$-MCT game on $P_2 \cp P_m$. 
	\paragraph{Case 1:} $2 \le m \le 5$.\\
	To show that $a_1(P_2 \cp P_m) \ge 2$, we will provide a winning strategy for Bob in the Alice-start $(1,1)$-MCT game. 
	Assume that Alice plays $(i,j)$. Then Bob replies by playing a vertex $(i',j')$ such that $i' \neq i$ and $j'\ne j$, and he can win the game with his next move. Thus $a_1(P_2 \cp P_m) \ge 2$. 
	
It remains to show that $a_1(P_2 \cp P_m) \le 2$. The assertion is clear for $m = 2$. Let $X_i= \{(1,2),(2,i-1) \}$, $i\in \{3,4,5\}$. Then $\max\{\Delta((P_2 \cp P_i) -X_i), |X_i|\} = 2$ holds for each $i\in \{3,4,5\}$. Hence Theorem~\ref{thm:threshold}(i) implies the assertion. 

	\paragraph{Case 2:} $ m \ge 6$.\\
It is straightforward to verify that by deleting any two vertices of $P_2\cp P_m$, the maximum degree of the rest of the graph remains $3$. Hence Theorem~\ref{thm:threshold}(i) yields $a_1(P_2 \cp P_m) \ge 3$. On the other hand, since  each vertex in $P_2 \cp P_m$ is incident with at most three cliques, Observation~\ref{Breakerwin}(iii) implies $a_1(P_2 \cp P_m) \le 3$. 
\qed 
	
	\begin{lemma}\label{lem:PnPmgamma}
		Assume that $m \ge n\ge 3$. Then $a_1(P_n \cp P_m) \le 3$ holds if and only if $\gamma(P_{n-2} \cp P_{m-2}) \le 3$.
	\end{lemma}
	
	\proof 
	The proof proceeds along the same lines as the corresponding inequality part of the proof of Lemma~\ref{lem:CnPm3}. 
	\qed

Obviously,  $\gamma(P_n) = \gamma(P_n\cp P_1) = \lceil\frac{n}{3}\rceil$. In~\cite{alanko-2011}, $\gamma(P_n \cp P_m)$ was calculated for $m,n \le 29$ by using a dynamic programming algorithm.  According to these results we may conclude the following:
	
	\begin{proposition}
		\label{prop:gamma-PnPm}	
		If $m \ge n\ge 1$, then $\gamma(P_n \cp P_m) \le 3$ if and only if one of the following conditions holds:
		\begin{itemize}
			\item $n=1$ and $1 \le m \le 9$;
			\item  $n=2$ and $2 \le m \le 5$;
			\item $n=m=3$.
		\end{itemize} 
	\end{proposition}	

\begin{lemma}\label{lem:PnPm2}
	Assume that $m \ge n\ge 3$. Then $a_1(P_n \cp P_m) = 2$ if and only if $n=3$ and $m\in \{3, 4\}$.
\end{lemma}

\proof 
If  $n=3$ and $m \ge 5$, or $ m \ge n \ge 4$, then by Theorem~\ref{thm:threshold}(i) we get $a_1(P_n \cp P_m) > 2$. Assume next that $n=3$ and $m\in \{3, 4\}$. Let $X = \{(2,2), (2,3)\}$. Then $\Delta((P_n \cp P_m) - X) = 2$, so $a_1(P_n \cp P_m) \le 2$ by Theorem~\ref{thm:threshold}(i). It remains to show that $a_1(P_n \cp P_m) \ge 2$ by providing a winning strategy for Bob in the Alice-start $(2,1)$-MCT game. After Alice plays $(i,j)$ in her first move, it is easy to see that there is $(i',j')$ which is $i' \neq i$ and $j' \neq j$. Then Bob replies by playing $(i',j')$ and wins the game in his next move. Hence $a_1(P_n \cp P_m) = 2$.	
\qed

For $m \ge n \ge 3$, according to Lemmas~\ref{lem:PnPmgamma} and \ref{lem:PnPm2}, and to Proposition~\ref{prop:gamma-PnPm}, we can conclude that $a_1(P_n \cp P_m) = 3$ if $m$ and $n$ satisfy one of the following conditions:
\begin{itemize}
	\item $n=3$ and $ 5 \ge m \ge 11$,
	\item $n=4$ and $ 4 \ge m \ge 7$,
	\item $n=m=5$.
\end{itemize}

To complete the proof of Theorem~\ref{thm:PnPm}, it remains to detect  when $m \ge n\ge 3$ and $a_1(C_n\cp P_m) > 3$ hold. By Observation~\ref{Breakerwin}(iii) we have  $a_1(P_n\cp P_m) \le 4$, hence in these cases $a_1(P_n\cp P_m) = 4$ holds. By Lemma~\ref{lem:PnPmgamma} and Proposition~\ref{prop:gamma-PnPm}, these are the following cases: 
\begin{itemize}
	\item $n = 3$ and $m \ge 12$, 
	\item $n = 4$ and $m \ge 8$,
	\item $n = 5$ and $m \ge 6$,
	\item $m \ge n \ge 6$.
\end{itemize}

\section*{Acknowledgments}

Sandi Klav\v{z}ar was supported by the Slovenian Research Agency (ARRS) under the grants P1-0297, J1-2452, and N1-0285. Csilla Bujt\'as acknowledges the support of the Slovenian Research Agency (ARRS) under the grant N1-0108.


\begin{thebibliography}{99}
		
		\bibitem{alanko-2011}
		S.~Alanko, S.~Crevals, A.~Isopoussu, P.~Ostergard, V.~Pettersson,
		Computing the domination number of grid graphs,
		Electron.\ J.\ Comb.\ 18 (2011) \#P141. 
		
		\bibitem{berge}
		C.~Berge,
		Hypergraphs. Combinatorics of Finite Sets,
		North Holland, 1989.
		
		\bibitem{bujtas-2022}
			Cs.~Bujt\'as, P.~Dokyeesun, Fast winning strategies for Staller in the Maker-Breaker domination game, http://arxiv.org/abs/2206.12812, 2021.
	
	\bibitem{campos-2013}
	C.~N.~Campos, S.~Dantas, C.~P.~De Mello, 
	Colouring clique-hypergraphs of circulant graphs,
	Graphs Comb.\ 29 (2013) 1713--1720.


		\bibitem{caro-2013}
		Y.~Caro, A.~Hansberg, New approach to the $k$-independence number of a graph, Electron.\ J.\ Combin.\ 20 (2013) \#P33.
		
		\bibitem{chellali-survey}
		M.~Chellali, O.~Favaron, A.~Hansberg, L.~Volkmann, $k$-Domination and $k$-independence in graphs: A survey, Graphs Combin.\ 28 (2012) 1--55.
		
	\bibitem{duchene-2020}
    E.~Duch\^{e}ne, V.~Gledel, A.~Parreau, G.~Renault,
    Maker-Breaker domination game,
    Discrete Math.\ 343 (2020) 111955.
    
		\bibitem{erdos-1973} 
		P.~Erd\H{o}s, J.L.~Selfridge, 
		On a combinatorial game, 
		J.\ Combin.\ Theory Ser.\ A 14 (1973) 298--301. 
		
	\bibitem{glazik-2022}
	 C.~Glazik, A.~Srivastav, 
	 A new bound for the Maker-Breaker triangle game,
	 European\ J.\ Combin.\ 104 (2022) 103536. 

	\bibitem{gledel-2020}
    V.~Gledel, M.~A.~Henning, V.~Ir\v si\v c, S.~Klav\v zar, 
    Maker-Breaker total domination game,
    Discrete Appl.\ Math.\ 282 (2020) 96--107.
	
		\bibitem{hefetz}
		D.~Hefetz, M.~Krivelevich, M.~Stojakovi\'c, T.~Szab\' o,
		Positional Games,
		Birkh\" auser/Springer, Basel, 2014.
		
		\bibitem{hik-2011}
		R.~Hammack, W.~Imrich, S.~Klav\v{z}ar,
		\textit{Handbook of Product Graphs, Second Edition},
		CRC Press, Boca Raton, FL, 2011.

	\bibitem{liang-2020}
    Z.~Liang, J.~Wu, E.~Shan, 
    List-coloring clique-hypergraphs of $K_5$-minor-free graphs strongly,
    Discrete Math.\ 343 (2020) 111777. 
	
		\bibitem{mao-2018}
		Y.~Mao, E.~Cheng, Z.~Wang, Z.~Guo, 
		The $k$-independence number of graph products,
		Art Discrete Appl.\ Math.\ 1 (2018) Paper No.\ P1.01. 
		
		\bibitem{mendes-2021}
		W.P.~Mendes, S.~Dantas, S.~Gravier, 
		The $(a, b)$-monochromatic transversal game on biclique-hypergraphs of powers of paths and of powers of cycles,
		Procedia Comp.\ Sci.\ 195 (2021) 181--189. 
		
		\bibitem{mendes-2022+}
		W.P.~Mendes, S.~Dantas, S.~Gravier, R.~Marinho,
		The $(a, b)$-monochromatic transversal game on clique-hypergraphs of powers of cycles,
		https://hal.archives-ouvertes.fr/hal-03015815, 2021.
		
		\bibitem{pavlic-2013}
		P.~Pavli\v{c}, J.~\v{Z}erovnik,
		A note on the domination number of the Cartesian products of paths and cycles,
		Kragujevac J.\  Math.\ 37 (2013) 275--285. 
		
	\bibitem{stojakovic-2021}
	 M.~Stojakovi\'{c}, N.~Trkulja, 
	 Hamiltonian Maker-Breaker games on small graphs,
	 Exp.\ Math.\ 30 (2021) 595--604.

		\bibitem{west-2001}
		D.B.~West,
		Introduction to Graph Theory, 2nd ed.,
		Prentice-Hall, NJ, 2001.
		
	\end{thebibliography}
\end{document}